\documentclass{amsart}
\allowdisplaybreaks
\usepackage{amscd}
\usepackage{amsfonts}
\usepackage{amsmath}
\usepackage{amssymb}
\usepackage{latexsym}
\usepackage{amsthm,color}
\usepackage[all]{xy}
\usepackage{epsfig}
\usepackage{graphicx}
\usepackage{color}
\usepackage{bm}
\usepackage{tikz}
\usetikzlibrary{cd}
\usetikzlibrary{calc}
\newtheorem{proposition}{Proposition}

\newtheorem{theorem}{Theorem}

\theoremstyle{definition}
\newtheorem{definition}{Definition}
\theoremstyle{definition}
\newtheorem{example}{Example}
\theoremstyle{definition}

\theoremstyle{definition}

\theoremstyle{remark}
\numberwithin{equation}{section}

\makeatletter
\@namedef{subjclassname@2020}{\textup{2020}Mathematics Subject Classification}
\makeatother

\begin{document}
	\title[
	Envelopes created by pseudo$\mbox{-}$circle families in the Minkowski plane]
	{Envelopes created by pseudo$\mbox{-}$circle families in the Minkowski plane}
	\author[Y.~Wang]{Yongqiao Wang
	}
	\address{
		School of Science, Dalian Maritime University, Dalian 116026, P.R. China
	}
	\email{wangyq@dlmu.edu.cn}
	\author[L.~Yang]{Lin Yang
	}
	\address{
		School of Science, Dalian Maritime University, Dalian 116026, P.R. China
	}
	\email{yanglin@dlmu.edu.cn}

    \author[Y. Chang]{Yuan Chang$^{*}$}
\address{School of Mathematics, Dongbei University of Finance and Economics, Dalian 116025, P.R. China}
\email{changy912@nenu.edu.cn}
\author[P. Li]{Pengcheng Li}
\address{Department of Mathematics, Harbin Institute of Technology, Weihai, 264209, P.R. China}
\email{lipc@hit.edu.cn}

	\begin{abstract}
	In this paper, we address the topic of envelopes created by pseudo$\mbox{-}$circle families in the Minkowski plane, which exhibit some different properties when compared with the Euclidean case. We provide solutions to all four basic problems associated with these envelopes, namely the existence problem, representation problem, problem on the number of envelopes, and problem on the relationships of definitions.	
	\end{abstract}


\date{}

\maketitle

\section{Introduction\label{section1}}
In classical differential geometry, regular submanifolds are widely studied. However, if a submanifold is not regular at some points, the classical ways studying the differential geometry of submanifold at those points do not apply. Fortunately, T. Fukunaga and M. Takahashi introduced a moving frame in the unit tangent bundle for studying frontals in \cite{fukunagatakahashi}. The main difference between a regular curve and a frontal is that the frontal may have singular points. In addition, by using this moving frame, they also investigate the relationship between Legendre curves and the arc-length parameter of Legendre immersions in the unit tangent bundle. In \cite{MT}, M. Takahashi defined an envelope for the $r$-parameter family of Legendre mappings and investigate their properties. As an application, they gave a condition that the projection of a singular solution of a first-order partial differential equation serves as an envelope.

This paper is a part of our research projects about the differential geometry of envelopes. The theory of envelopes plays a crucial role in various areas such as differential geometry, differential equations, geometric optics and so on.
In \cite{nishimura}, T. Nishimura studied the envelopes for a given hyperplane family and developed a comprehensive theory for constructing the envelope of the hyperplane family. However, the family of circles on a plane should not be disregarded  because
their envelopes have already been extensively applied, such as in the field of seismic survey, see \cite{brucegiblin}.
In \cite{Wt}, the first author and T. Nishimura studied the envelope created by a family of circles in the Euclidean plane and successfully addressed some key questions related to this envelope. These problems include the necessary and sufficient conditions for a family of circles to create an envelope, the parametric representation of the envelope, and the relationship between the number of envelopes and the central locus and radius function of the circle family.

On the other hand, the properties of submanifolds in the Minkowski space differ from those in Euclidean space, due to the deviation of the inner product from that of Euclidean space, exhibiting degeneracy and negativity. In \cite{AF}, A. Saloom and F. Tari studied the evolute of spacelike curves and timelike curves in the Minkowski plane. Additionally, they studied the generic behaviour of caustic of a smooth curve at the lightlike points.
In this paper, we investigate the envelope created by a pseudo$\mbox{-}$circle family in the Minkowski plane. As it to be expected, the envelope of the pseudo$\mbox{-}$circle family exhibits some different properties when compared with the Euclidean case. For instance, the pseudo$\mbox{-}$circle family may not create an envelope if the radius function is constant, see Propositions \ref{1.13} and \ref{1.133} in Section \ref{section2}.

Throughout this paper, unless otherwise specified, all manifolds, functions and mappings are assumed to be class $C^{\infty}$.

\section{Preliminaries and main results\label{section2}}
It is known that Minkowski plane $\mathbb{R}^2_1$ is the plane $\mathbb{R}^2$ with a pseudo-scalar product $\langle\bm{x},\bm{y}\rangle=-x_1y_1+x_2y_2$, where $\bm{x}=(x_1,x_2)$ and $\bm{y}=(y_1,y_2)$. By definition, this product
classifies the vectors in $\mathbb{R}^2_1$ as follows: $\bm{x}$ is said to be spacelike, lightlike, and timelike if $\langle\bm{x},\bm{x}\rangle>0,=0,<0,$ respectively.
Let $P$ be a point in the Minkowski plane $\mathbb{R}^2_1$ and let $r>0$, we consider following pseudo$\mbox{-}$circles in $\mathbb{R}^2_1$ with center at $P$ and radius $r$.
$$S^1_1(P, r)=\{\bm{x}\in\mathbb{R}^2_1|\langle \bm{x}-P, \bm{x}-P\rangle=r^2\},$$
$$H^1(P, -r)=\{\bm{x}\in\mathbb{R}^2_1|\langle \bm{x}-P, \bm{x}-P\rangle=-r^2\}.$$
We denote $S^1_1(0, 1)$ and $H^1(0, -1)$ as $S^1_1$ and $H^1$, respectively. We call $\bm{x}$ a unit vetor if $\bm{x}\in S^1_1~ (\rm{ or }~H^1)$. Let $\bm{a}:I\rightarrow\mathbb{R}^2_1$ be a curve in the Minkowski plane and $r:I\rightarrow\mathbb{R}^+
$ a positive function, where $\mathbb{R}^+$ is the set of all positive real numbers. As a natural generalization,  we define pseudo$\mbox{-}$circle family $C_{(\bm{a}(t), \pm r(t))}$ center at $a(t)$ with radius $r(t)$ as follows:
$$C_{(\bm{a}(t), \pm r(t))}=\{\bm{x}\in\mathbb{R}^2_1|\langle \bm{x}-\bm{a}(t), \bm{x}-\bm{a}(t)\rangle=\pm r(t)^2\}.$$

We say that $\bm{a}:I\rightarrow\mathbb{R}_1^2$ is a spacelike(\rm{or }~ timelike) frontal, if there exists a smooth map $\bm{\nu}:I\rightarrow H^1(\rm{ or }~S^1_1)$ such that the pair  $(\bm{a}, \bm{\nu}):I\times H^1(\rm{ or }~S^1_1)\rightarrow \mathbb{R}_1^2$ is a Legendre curve. It follows that the derivative of $\bm{a}(t)$ with respect to $t$, denoted as $\frac{d \bm{a}}{d t}(t)$, is pseudo$\mbox{-}$orthogonal to $\bm{\nu}(t)$ which is represented as $\langle\frac{d \bm{a}}{d t}(t), \bm{\nu}(t)\rangle=0$. The mapping $\bm{\nu}$ is commonly known as the Gauss mapping of frontal $\bm{a}$. There exists a unit vector $\bm{\mu}(t)$ such that $\langle\bm{\mu}(t),\bm{\nu}(t)\rangle=0$. Then,
the pair $\{\bm{\mu}, \bm{\nu}\}$ is a moving frame along frontal $\bm{a}$. We define the sign of a non$\mbox{-}$lightlike vector $\bm{x}$ as $\epsilon_{\bm{x}}=1$ if $\bm{x}$ is a spacelike vector and $\epsilon_{\bm{x}}=-1$ if $\bm{x}$ is a timelike vector, respectively. We have the following Frenet formula:
\begin{equation*}
	\frac{d }{d t}\left(\begin{array}{c} \bm{\nu}(t) \\ \bm{\mu}(t)\end{array}\right)
	=\left(\begin{array}{c c}0 & \ell(t) \\\ell(t) & 0\end{array}\right)
	\left(\begin{array}{c}\bm{\nu}(t) \\ \bm{\mu}(t)\end{array}\right),~~
	\frac{d \bm{a}}{d t}(t)=\beta(t)\bm{\mu}(t),
\end{equation*}
where  $\ell(t)=\epsilon_{\bm{\mu}}\langle\bm{\mu}(t),\frac{d \bm{\nu}}{d t}(t)\rangle$, $\beta(t)=\epsilon_{\bm{\mu}}\langle\bm{\mu}(t),\frac{d \bm{a}}{d t}(t)\rangle$. The pair $(\ell,\beta)$ is called the curvature of the frontal $\bm{a}$ (for details, see \cite{GS}).

By definition, a frontal is a mapping that provides a solution to the first-order linear differential equation defined by the Gauss mapping $\bm{\nu}$. Therefore, the set of frontals with a fixed mapping $\bm{\nu}$, characterized by a given Gauss mapping $\bm{\nu}:I\rightarrow H^1(\rm{ or }~S^1_1)$, forms a linear space. Then, we can infer that frontal is the solution of a first-order linear differential equation given by the mapping $\bm{\nu}$.
 In this paper, we assume that the curve $\bm{a}(t)$ for a pseudo$\mbox{-}$circle family $C_{(\bm{a}(t), \pm r(t))}$ is a frontal.
 Then, we have following definition for the envelope created by $C_{(\bm{a}(t), \pm r(t))}$.

\begin{definition}\label{1.15}
	Let $C_{(\bm{a}(t), \pm r(t))}$ be a pseudo$\mbox{-}$circle family in the Minkowski plane. A smooth mapping  $f:I\rightarrow\mathbb{R}_1^2$ is called an envelope of $C_{(\bm{a}(t), \pm r(t))}$ if $f(t)$ satisfies following two conditions:\\
	$(1)$ $f(t)\in C_{(\bm{a}(t), \pm r(t))},$\\
	$(2)$ $\langle \frac{d f}{d t}(t),f(t)-\bm{a}(t)\rangle=0$.
\end{definition}
By the above definition, an envelope created by a pseudo$\mbox{-}$circle family is a mapping that provides a solution to a first-order differential equation subject to a constraint condition.
The following Example \ref{1.111} shows that the well-known method to calculate the envelope created by a pseudo$\mbox{-}$circle family.
\begin{example}\label{1.111}
	Let $\bm{a}:I\rightarrow\mathbb{R}_1^2$ be a timelike frontal defined by $\bm{a}(t)=(t^3,\sqrt{1+t^6})$, where $t\in\mathbb{R}$. We calculate that $\bm{\nu}(t)=(t^3,\sqrt{1+t^6})$ and $\bm{\mu}(t)=(\sqrt{1+t^6},t^3)$. Let $r:\mathbb{R}\rightarrow\mathbb{R}^+$ be a positive function defined by $r(t)=1$. Then, it appears that the pseudo$\mbox{-}$circle family $C_{(\bm{a}(t), r(t))}$ creates envelopes. Thus, by direct calculation, the envelopes of $C_{(\bm{a}(t), r(t))}$ are obtained as
	\begin{align*}
		\mathcal{D} = &
		\left\{(x, y)\in \mathbb{R}^2_1\, \left|\,
		\exists t \mbox{ such that }F(x,y,t)=\frac{\partial F}{\partial t}(x, y, t)=0
		\right.\right\} \\
		 = &
		\left\{(x, y)\in \mathbb{R}^2_1\, \left|\,
		\exists t \mbox{ such that }
		-\left(x-t^3\right)^2+\left(y-\sqrt{1+t^6}\right)^2-1=0, \,
		6t^2\left(x-t^3\right)-\frac{6t^5}{\sqrt{1+t^6}}\left(y\right.\right.\right. \\
		&
		\left.\left.\left.-\sqrt{1+t^6}\right)=0
		\right.\right\} \\
		 = &
		\left\{(x, y)\in \mathbb{R}^2_1\, \left|\,
		\exists t \mbox{ such that }
		-\left(x-t^3\right)^2+\left(y-\sqrt{1+t^6}\right)^2-1=0, \,
		t^2\left(\left(x-t^3\right)-\frac{t^3}{\sqrt{1+t^6}}\left(y\right.\right.\right.\right. \\
		 &
		\left.\left.\left.\left.-\sqrt{1+t^6}\right)\right)=0
		\right.\right\} \\
		= &
		\left\{(x, y)\in \mathbb{R}^2_1\, \left|\,
		-x^2+y^2=1\right.\right\} \bigcup
		\left\{(x, y)\in \mathbb{R}^2\, \left|\,
		\left(x-t^3\right)^2-\left(y-\sqrt{1+t^6}\right)^2+1=0, \,
		\right.\right. \\
		  &
		\left.\left.x=t^3+\frac{t^3}{\sqrt{1+t^6}}\left(y-\sqrt{1+t^6}\right)
		\right.\right\} \\
		= &
		\left\{(x, y)\in \mathbb{R}^2_1\, \left|\,
		-x^2+y^2=1\right.\right\} \bigcup
		\left\{(x, y)\in \mathbb{R}^2\, \left|\,
		\left(\frac{t^3}{\sqrt{1+t^6}}\left(y-\sqrt{1+t^6}\right)\right)^2-\left(y-\sqrt{1+t^6}\right)^2=1, \,
		\right.\right. \\
		  &
		\left.\left.
		x=\frac{yt^3}{\sqrt{1+t^6}}
		\right.\right\}
		\\
		 = &
		\left\{(x, y)\in \mathbb{R}^2_1\, \left|\,
		-x^2+y^2=1\right.\right\} \bigcup
		\left\{\left.
		\left(t^3\pm t^3,\sqrt{1+t^6}\pm \sqrt{1+t^6}\right)
		\in \mathbb{R}^2_1\, \right|\, t\in \mathbb{R}
		\right\}.
	\end{align*}
	\begin{figure}[h]
		\begin{center}
			\includegraphics[width=8cm]
			{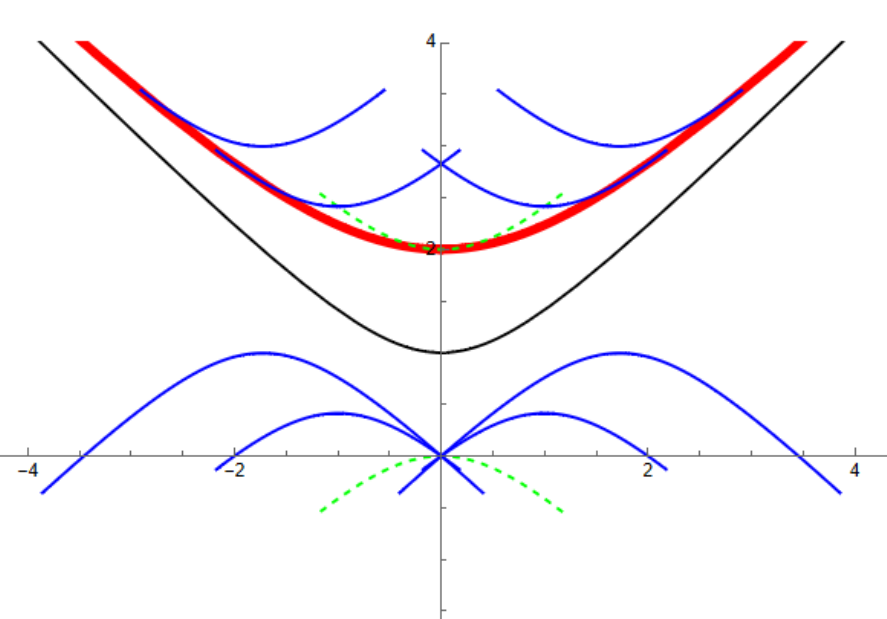}
			\caption{The pseudo$\mbox{-}$circle family $C_{(\bm{a}(t), r(t))}$
				and the candidates of its envelopes.
			}
			\label{figure_example1}
		\end{center}
	\end{figure}
\end{example}
The set $\mathcal{D}$ obtained by the classical method is called as $\mathcal{D}$ \textit{envelope} in Definition \ref{4.1} of Section \ref{section4}. In Example \ref{3.1} of Section \ref{section3}, it is revealed that the $\mathcal{D}$ \textit{envelope} calculated here is indeed larger than the set of envelopes created by $C_{(\bm{a}(t), r(t))}$, making the inclusion of the unit pseudo$\mbox{-}$circle $\left\{(x, y)\in \mathbb{R}^2_1\, \left|\,-x^2+y^2=1\right.\right\}$ redundant, (see Figure \ref{figure_example1}, the dashed curves are redundant). The classical method
to acquire the envelopes does not work well in this case. In order to accurately calculate the envelopes of the pseudo$\mbox{-}$circle family, we have the following key definition.
\begin{definition}\label{1.1111}
	Let $(\bm{a}, \bm{\nu}):I\times  H^1(\rm{ or }~S^1_1)\rightarrow\mathbb{R}_1^2$ be a spacelike(\rm{or }~  timelike) frontal and $r:I\rightarrow\mathbb{R}^+$ a positive function. Then, the family $C_{(\bm{a}(t), r(t))}$ is said to be creative if there exist smooth mapping $\tilde{\bm{\nu}}:I\rightarrow  S^1_1$ such that
	$$\frac{d r}{dt}(t)+\beta(t)\langle\tilde{\bm{\nu}}(t),\bm{\mu}(t)\rangle=0.$$
	The family $C_{(\bm{a}(t), -r(t))}$ is said to be creative if there exist a smooth mapping $\tilde{\bm{\nu}}:I\rightarrow  H^1$ such that
	$$\frac{d r}{dt}(t)-\beta(t)\langle\tilde{\bm{\nu}}(t),\bm{\mu}(t)\rangle=0.$$

If $\epsilon_{\tilde{\bm{\nu}}}\cdot\epsilon_{\bm{\mu}}=1$, the creative condition is equivalent to the existence of a function $\theta:I\rightarrow\mathbb{R}$ satisfying the following identities for any $t\in I$, $$\frac{d r}{dt}(t)+\beta(t)\cosh\theta(t)=0$$ or
$$\frac{d r}{dt}(t)-\beta(t)\cosh\theta(t)=0.$$

If $\epsilon_{\tilde{\bm{\nu}}}\cdot\epsilon_{\bm{\mu}}=-1$, the creative condition is equivalent to the existence of a function $\theta:I\rightarrow\mathbb{R}$ satisfying the following identity for any $t\in I$, $$\frac{d r}{dt}(t)-\beta(t)\sinh\theta(t)=0.$$
\end{definition}
Next, we give the necessary and sufficient conditions for a family of pseudo$\mbox{-}$circles to create an envelope. We have following theorem.
\begin{theorem}\label{1.14}
	Let $\bm{a}:I\rightarrow\mathbb{R}_1^2$ be a frontal with the Gauss mapping $\bm{\nu}:I\rightarrow H^1(\rm{ or }~S^1_1)$ and let $r:I\rightarrow\mathbb{R}^+$ be a positive function. Then, the pseudo$\mbox{-}$circle family $C_{(\bm{a}(t), \pm r(t))}$ creates an envelope if and only if $C_{(\bm{a}(t), \pm r(t))}$ is creative.
\end{theorem}
\begin{proof}
	We suppose that $C_{(\bm{a}(t),\pm r(t))}$ is creative. By definition, there exists a mapping $\tilde{\bm{\nu}}:I\rightarrow S^1_1(\rm{or }~ H^1)$ such that
	$\frac{d r}{dt}(t)=\mp\beta(t)\langle\tilde{\bm{\nu}}(t),\bm{\mu}(t)\rangle$ holds for any $t\in I$. Let $f:I\rightarrow\mathbb{R}_1^2$ be a curve defined by
	$$f(t)=\bm{a}(t)+r(t)\tilde{\bm{\nu}}(t).$$
	We have $\langle f(t)-\bm{a}(t),f(t)-\bm{a}(t)\rangle=r^2(t)\langle\tilde{\bm{\nu}}(t),\tilde{\bm{\nu}}(t)\rangle$, where $\langle\tilde{\bm{\nu}}(t),\tilde{\bm{\nu}}(t)\rangle=\pm1$.
	It means that $f(t)\in C_{(\bm{a}(t),\pm r(t))}$. Morever, since
	$$\frac{d f}{dt}(t)=\frac{d \bm{a}}{dt}(t)+\frac{d r}{dt}(t)\tilde{\bm{\nu}}(t)+r(t)\frac{d \tilde{\bm{\nu}}}{dt}(t),$$
	we have
	\begin{align*}
		&\big\langle \frac{d f}{dt}(t),f(t)-\bm{a}(t)\big\rangle\\
=&\big\langle \frac{d \bm{a}}{dt}(t)+\frac{d r}{dt}(t)\tilde{\bm{\nu}}(t)+r(t)\frac{d \tilde{\bm{\nu}}}{dt}(t),r(t)\tilde{\bm{\nu}}(t)\big\rangle\\
=&\big\langle \frac{d \bm{a}}{dt}(t),r(t)\tilde{\bm{\nu}}(t)\big\rangle+\big\langle \frac{d r}{dt}(t)\tilde{\bm{\nu}}(t),r(t)\tilde{\bm{\nu}}(t)\big\rangle\\
=&\langle \beta(t)\bm{\mu}(t),r(t)\tilde{\bm{\nu}}(t)\rangle\pm \frac{d r}{dt}(t)r(t)\\
=&r(t)\beta(t)\langle \bm{\mu}(t),\tilde{\bm{\nu}}(t)\rangle\pm(\mp\beta(t)\langle\tilde{\bm{\nu}}(t),\bm{\mu}(t)\rangle)r(t)\\
=&0.		
	\end{align*}
	Therefore, $f(t)$ is an envelope created by the pseudo$\mbox{-}$circle family $C_{(\bm{a}(t),\pm r(t))}$.
	
	Conversely, suppose that the pseudo$\mbox{-}$circle family $C_{(\bm{a}(t),\pm r(t))}$ creates an envelope $f:I\rightarrow\mathbb{R}_1^2$. Then, it implies that $f(t)\in C_{(\bm{a}(t), \pm r(t))}$ and $\langle \frac{d f}{dt}(t),f(t)-\bm{a}(t)\rangle=0$. The condition $f(t)\in C_{(\bm{a}(t), \pm r(t))}$ follows that there exists a mapping $\tilde{\bm{\nu}}:I\rightarrow S^1_1(\rm{or }~ H^1)$ such that the following equality holds for any $t\in I$,
	$$f(t)=\bm{a}(t)+r(t)\tilde{\bm{\nu}}(t).$$
	Since
	$$\frac{d f}{dt}(t)=\frac{d \bm{a}}{dt}(t)+\frac{d r}{dt}(t)\tilde{\bm{\nu}}(t)+r(t)\frac{d \tilde{\bm{\nu}}}{dt}(t),$$
	we have
	\begin{align*}
		0&=\langle \frac{d f}{dt}(t),f(t)-\bm{a}(t)\rangle\\
		&=\langle \frac{d \bm{a}}{dt}(t)+\frac{d r}{dt}(t)\tilde{\bm{\nu}}(t)+r(t)\frac{d \tilde{\bm{\nu}}}{dt}(t),r(t)\tilde{\bm{\nu}}(t)\rangle\\
		&=\langle \beta(t)\bm{\mu}(t),r(t)\tilde{\bm{\nu}}(t)\rangle\pm \frac{d r}{dt}(t)r(t)\\
		&=r(t)\big(\langle \beta(t)\bm{\mu}(t),\tilde{\bm{\nu}}(t)\rangle\pm \frac{d r}{dt}(t)\big).		
	\end{align*}
	Because $r:I\rightarrow\mathbb{R}^+$ is a positive function for any $t\in I$, we have
	$$\frac{d r}{dt}(t)\pm\langle \beta(t)\bm{\mu}(t),\tilde{\bm{\nu}}(t)\rangle=0.$$
	Therefore, the pseudo$\mbox{-}$circle family $C_{(\bm{a}(t),\pm r(t))}$ is creative.
\end{proof}
In order to obtain the parametric representation of the envelope created by a creative family $C_{(\bm{a}(t),\pm r(t))}$, we have the following theorem.
\begin{theorem}\label{1.11}
	Let $\bm{a}:I\rightarrow\mathbb{R}_1^2$ be a frontal with the Gauss mapping $\bm{\nu}:I\rightarrow H^1(\rm{or }~ S^1_1)$ and let $r:I\rightarrow\mathbb{R}^+$ be a positive function. Suppose that the pseudo$\mbox{-}$circle family $C_{(\bm{a}(t),\pm r(t))}$ creates an envelope $f:I\rightarrow\mathbb{R}_1^2$. Then, the following two hold.\\
	$(1)$ If  $\epsilon_{\tilde{\bm{\nu}}}\cdot\epsilon_{\bm{\mu}}=1$, the envelope $f(t)$ is represented by
	$$f(t)=\bm{a}(t)+r(t)\tilde{\bm{\nu}}(t),$$
	where $\tilde{\bm{\nu}}(t)=\cosh\theta(t)\bm{\mu}(t)\pm\sinh\theta(t)\bm{\nu}(t)$ or $\tilde{\bm{\nu}}(t)=-\cosh\theta(t)\bm{\mu}(t)\pm\sinh\theta(t)\bm{\nu}(t)$, and $\theta(t)$ is the function in Definition \ref{1.1111}.\\
	$(2)$ If $\epsilon_{\tilde{\bm{\nu}}}\cdot\epsilon_{\bm{\mu}}=-1$, the envelope $f(t)$ is represented by
	$$f(t)=\bm{a}(t)+r(t)\tilde{\bm{\nu}}(t),$$
	where $\tilde{\bm{\nu}}(t)=\sinh\theta(t)\bm{\mu}(t)\pm\cosh\theta(t)\bm{\nu}(t)$ and $\theta(t)$ is the function in Definition \ref{1.1111}.
\end{theorem}
\begin{proof}
	The proof of the Theorem \ref{1.14} proves the Theorem \ref{1.11} as well.
\end{proof}
According to Definition \ref{1.1111} and Theorem \ref{1.14}, we have the following two propositions.
\begin{proposition}\label{1.13}
	Suppose that $\bm{a}:I\rightarrow\mathbb{R}_1^2$ is a spacelike frontal with the Gauss mapping $\bm{\nu}:I\rightarrow H^1$ and $r:I\rightarrow\mathbb{R}^+$ is a positive constant. If $\bm{a}(t)$ is not a point, then the pseudo$\mbox{-}$circle family $C_{(\bm{a}(t), r(t))}$ does not create an envelope.
\end{proposition}
\begin{proof}
	Since $r(t)$ is a positive constant, we have $\frac{d r}{dt}(t)=0$. By Definition \ref{1.1111}, the creative condition is equivalent to the existence of a function $\theta:I\rightarrow\mathbb{R}$ satisfying the following identities for any $t\in I$, $$\frac{d r}{dt}(t)+\beta(t)\cosh\theta(t)=0,$$ or
	$$\frac{d r}{dt}(t)-\beta(t)\cosh\theta(t)=0.$$
	Because $\bm{a}(t)$ is not a point, then $\beta(t)$ is not always equal to 0. It follows from $\cosh\theta(t)\geq1$ for any $\theta(t)$ that
	the above equations do not hold. Then the pseudo$\mbox{-}$circle family $C_{(\bm{a}(t), r(t))}$ is not creative. By Theorem \ref{1.14}, $C_{(\bm{a}(t), r(t))}$ does not create envelopes.
\end{proof}
\begin{proposition}\label{1.133}
	Suppose that $\bm{a}:I\rightarrow\mathbb{R}_1^2$ is a timelike frontal with the Gauss mapping $\bm{\nu}:I\rightarrow S^1_1$ and $r:I\rightarrow\mathbb{R}^+$ is a positive constant, If $\bm{a}(t)$ is not a point, then the pseudo$\mbox{-}$circle family $C_{(\bm{a}(t), -r(t))}$ does not create an envelope.
\end{proposition}
\begin{proof}
The proof is similar to the proof of Proposition \ref{1.13}.	
\end{proof}
Now we want to focus on $\frac{d r}{dt}(t)$ and $\beta(t)$ to study the number of envelopes created by a creative family $C_{(\bm{a}(t), \pm r(t))}$. We have the following theorem.
\begin{theorem}\label{1.12}
	Let $a:I\rightarrow\mathbb{R}_1^2$ be a frontal with the Gauss mapping $\bm{\nu}:I\rightarrow H^1(\rm{or }~ S^1_1)$ and $r:I\rightarrow\mathbb{R}^+$ a positive function. Suppose that the pseudo$\mbox{-}$circle family $C_{(\bm{a}(t), \pm r(t))}$ create an envelope. Then, the number of envelopes created by $C_{(\bm{a}(t), \pm r(t))}$ can be described as follows:\\
	$(1)$ Suppose that $\epsilon_{\tilde{\bm{\nu}}}\cdot\epsilon_{\bm{\mu}}=1$.
	
	$(1\mbox{-}i)$ The pseudo$\mbox{-}$circle family $C_{(\bm{a}(t), \pm r(t))}$ creates a unique envelope if and only if the set consisting of $t\in I$ satisfying $\beta(t)\neq0$ and $|\frac{d r}{dt}(t)|=|\beta(t)|$ is dense in $I$.
	
	$(1\mbox{-}ii)$ There are exactly two envelopes created by $C_{(\bm{a}(t), \pm r(t))}$ if and only if the set of $t\in I$ satisfying $\beta(t)\neq0$ is dense in $I$ and there exists at least one $t_0\in I$ such that $|\frac{d r}{dt}(t_0)|>|\beta(t_0)|$.
	
	$(1\mbox{-}\infty)$ There are uncountably many distinct envelopes created by $C_{(\bm{a}(t), \pm r(t))}$ if and only if the set of $t\in I$ satisfying $\beta(t)\neq0$ is not dense in $I$.\\
	$(2)$ Suppose that $\epsilon_{\tilde{\bm{\nu}}}\cdot\epsilon_{\bm{\mu}}=-1$.

$(2\mbox{-}i)$ It is impossible for the pseudo$\mbox{-}$circle family $C_{(\bm{a}(t), \pm r(t))}$ creating a unique envelope.

$(2\mbox{-}ii)$ There are exactly two envelopes created by $C_{(\bm{a}(t), \pm r(t))}$ if and only if the set of $t\in I$ satisfying $\beta(t)\neq0$ is dense in $I$.

$(2\mbox{-}\infty)$ There are uncountably many distinct envelopes created by $C_{(\bm{a}(t), \pm r(t))}$ if and only if the set of $t\in I$ satisfying $\beta(t)\neq0$ is not dense in $I$.\\
\end{theorem}
\begin{proof}
	$(1\mbox{-}i)$ Suppose that the pseudo$\mbox{-}$circle family $C_{(\bm{a}(t), \pm r(t))}$ create a unique envelope. By Theorem \ref{1.11}, the set $\{f:I\rightarrow\mathbb{R}_1^2~|~f(t)=\bm{a}(t)+r(t)\tilde{\bm{\nu}}(t)~~(t\in I)\}$  consists of only one mapping. Hence, for any $t\in I$ the unit vector $\tilde{\bm{\nu}}(t)$ satisfying
	$$\frac{d r}{dt}(t)+\epsilon_{\tilde{\bm{\nu}}}\beta(t)\langle\tilde{\bm{\nu}}(t),\bm{\mu}(t)\rangle=0$$
	must be uniquely determined. Since $\epsilon_{\tilde{\bm{\nu}}}\cdot\epsilon_{\bm{\mu}}=1$, by Definition \ref{1.1111}, we have $|\langle\tilde{\bm{\nu}}(t),\bm{\mu}(t)\rangle|=\cosh\theta(t)\geq1$. Hence, under considering continuity of two functions $\frac{d r}{dt}(t)$ and $\beta(t)$, it means that the set consisting of $t\in I$ satisfying $|\frac{d r}{dt}(t)|=|\beta(t)|\neq0$ is dense in $I$.
	
	Conversely, suppose that the set consisting of $t\in I$ satisfying $|\frac{d r}{dt}(t)|=|\beta(t)|\neq0$ is dense in $I$. Then, under considering continuity of the function $t\mapsto\langle\tilde{\bm{\nu}}(t),\bm{\mu}(t)\rangle$, it follows that for any $t\in I$, $\langle\tilde{\bm{\nu}}(t),\bm{\mu}(t)\rangle=1$(resp.,~ $\langle\tilde{\bm{\nu}}(t),\bm{\mu}(t)\rangle=-1$) if $\frac{d r}{dt}(t)=-\beta(t)$ and $\bm{\mu}(t)\in S^1_1$ (resp.,~$\bm{\mu}(t)\in H^1$), or
	$\langle\tilde{\bm{\nu}}(t),\bm{\mu}(t)\rangle=-1$(resp.,~ $\langle\tilde{\bm{\nu}}(t),\bm{\mu}(t)\rangle=1$) if $\frac{d r}{dt}(t)=\beta(t)$ and $\bm{\mu}(t)\in S^1_1$ (resp.,~$\bm{\mu}(t)\in H^1$).
	 Thus, the mapping $\tilde{\bm{\nu}}$ is unique. By Theorem \ref{1.11}, the envelope $f:I\rightarrow\mathbb{R}_1^2$ must be uniquely determined.\\
	\\
	$(1\mbox{-}ii)$ Suppose that there are exactly two distinct envelopes created by $C_{(\bm{a}(t), \pm r(t))}$. Then, by the equality $\frac{d r}{dt}(t)+\epsilon_{\tilde{\bm{\nu}}}\beta(t)\langle\tilde{\bm{\nu}}(t),\bm{\mu}(t)\rangle=0$, the set consisting of $t\in I$ satisfying $\beta(t)\neq0$ must be dense in $I$. Suppose moreover that the equality $|\frac{d r}{dt}(t)|=|\beta(t)|\neq0$ holds for any $t\in I$. Then it follows that the set consisting of $t\in I$ satisfying $|\frac{d r}{dt}(t)|=|\beta(t)|\neq0$ is dense in $I$. Then, by the assertion $(1\mbox{-}i)$, the given pseudo$\mbox{-}$circle family must create a unique envelope. This contradicts the assumption that there are exactly two distinct envelope. Hence, there must exist at least one $t_0\in I$ such that the strict inequality $|\frac{d r}{dt}(t_0)|>|\beta(t_0)|$ holds.
	
	Conversely, suppose that the set of $t\in I$ satisfying $\beta(t)\neq0$ is dense in $I$ and there exists at least one $t_0\in I$ such that the strict inequality $|\frac{d r}{dt}(t_0)|>|\beta(t_0)|$ holds. Then, it follows that there must exist an open intervl $\tilde{I}$ in $I$ such that the absolute value $|\langle\tilde{\bm{\nu}}(t),\bm{\mu}(t)\rangle|=\cosh\theta(t)>1$ for any $t\in \tilde{I}$. Thus, it follows $\theta(t)\neq-\theta(t)$ for any $t\in \tilde{I}$. Hence, for any $t\in \tilde{I}$, there exist exactly two distinct unit vectors $\tilde{\bm{\nu}}_+(t)$ and $\tilde{\bm{\nu}}_-(t)$ correspeonding $|\langle\tilde{\bm{\nu}}_+(t),\bm{\mu}(t)\rangle|=\cosh\theta(t)$ and $|\tilde{\bm{\nu}}_-(t),\bm{\mu}(t)\rangle|=\cosh(-\theta(t))$, respectively. Therefore, by the assertion of Theorem \ref{1.11}, the pseudo$\mbox{-}$circle family must create exactly two distinct envelopes.\\
	\\
	$(1\mbox{-}\infty)$ Suppose that there are uncountably many distinct  envelopes created by $C_{(\bm{a}(t), \pm r(t))}$. Suppose moreover that the  set of $t\in I$ such that $\beta(t)\neq0$ is dense in $I$. Then, from $(1\mbox{-}i)$ and $(1\mbox{-}ii)$, it follows that the pseudo$\text{-}$circle family $C_{(\bm{a}(t), \pm r(t))}$ must create a unique envelope or two distinct envelopes. This contradicts the assumption that there are uncountably many distinct  envelopes created by $C_{(\bm{a}(t), \pm r(t))}$. Hence, the set of $t\in I$ such that $\beta(t)\neq0$ is never dense in $I$.
	
	Conversely, suppose that the set of $t\in I$ such that $\beta(t)\neq0$ is not dense in $I$. This assumption implies that there exists an open interval $\tilde{I}$ in $I$ such that $\beta(t)=0$ for any $t\in\tilde{I}$. On the other hand, since $C_{(\bm{a}(t), \pm r(t))}$ creates an envelope $f_0$, the equality
	$$\frac{d r}{dt}(t)+\epsilon_{\tilde{\bm{\nu}}}\beta(t)\langle\tilde{\bm{\nu}}(t),\bm{\mu}(t)\rangle=0$$
	holds for any $t\in I$. Thus, there are no restrictions for the value $\langle\tilde{\bm{\nu}}(t),\bm{\mu}(t)\rangle$ for any $t\in \tilde{I}$. Take one point $t_0$ of $\tilde{I}$ and denote the $\tilde{\bm{\nu}}(t)$ for the envelope $f_0$ by $\tilde{\bm{\nu}}_0(t)$. Then, by using the standard technique on bump functions, we may construct uncountably many distinct unit vectors $\tilde{\bm{\nu}}_a(t)\in S^1_1(\rm{or }~ H^1)$ that satisfy the following four conditions, where $a\in A$ and $A$ is a set consisting uncountably many elements such that $0\notin A$.
	
	$(a)$ For any $t\in I$ and any $a\in A$, the equality $\frac{d r}{dt}(t)=-\epsilon_{\tilde{\bm{\nu}}_a}\beta(t)\langle\tilde{\bm{\nu}}_a(t),\bm{\mu}(t)\rangle$ holds.
	
	$(b)$ For any $t\in I-\tilde{I}$ and any $a\in A$, the equality $\tilde{\bm{\nu}}_a(t)=\tilde{\bm{\nu}}_0(t)$ holds.
	
	$(c)$ For any $t\in\tilde{I}$ and $a\in A$, $\tilde{\bm{\nu}}_a(t_0)\neq\tilde{\bm{\nu}}_0(t_0)$ holds.
	
	$(d)$ For any $a_1, a_2\in A$ and $a_1\neq a_2$, $\tilde{\bm{\nu}}_{a_1}(t_0)\neq\tilde{\bm{\nu}}_{a_2}(t_0)$ holds.\\
	Therefore, by Theorem \ref{1.11}, the pseudo$\text{-}$circle $C_{(\bm{a}(t), \pm r(t))}$ creates uncountably many distinct envelopes.\\
	\\	
	$(2\mbox{-}i)$ Suppose that the pseudo$\mbox{-}$circle family $C_{(\bm{a}(t), \pm r(t))}$ creates an envelope. Since the creative condition is equivalent to the existence of a function $\theta:I\rightarrow\mathbb{R}$ satisfying $\frac{d r}{dt}(t)-\beta(t)\sinh\theta(t)=0$ for any $t\in I$, we divide two cases: $(A)$ $\beta(t)\neq0$ is dense in $I$, $(B)$ $\beta(t)\neq0$ is not dense in $I$.
	
	First, we consider the case $(A)$. Since $\beta(t)\neq0$ is dense in $I$, by Definition \ref{1.1111} the unit vector $\tilde{\bm{\nu}}$ satisfying $$\frac{d r}{dt}(t)+\epsilon_{\tilde{\bm{\nu}}}\beta(t)\langle\tilde{\bm{\nu}}(t),\bm{\mu}(t)\rangle=0$$ for any $t\in I$. It follows $\langle\tilde{\bm{\nu}}(t),\bm{\mu}(t)\rangle=\epsilon_{\bm{\mu}}\sinh\theta(t)$ for any $t\in I$. Since $\cosh\theta(t)\neq-\cosh\theta(t)$, there exist exactly two distinct unit vectors $\tilde{\bm{\nu}}_+(t)=\sinh\theta(t)\bm{\mu}(t)+\cosh\theta(t)\bm{\nu}(t)$ and $\tilde{\bm{\nu}}_-(t)=\sinh\theta(t)\bm{\mu}(t)-\cosh\theta(t)\bm{\nu}(t)$
	such that $\langle\tilde{\bm{\nu}}_+(t),\bm{\mu}(t)\rangle=\langle\tilde{\bm{\nu}}_-(t),\bm{\mu}(t)\rangle=\epsilon_{\bm{\mu}}\sinh\theta(t)$. Therefore, by the assertion of Theorem \ref{1.11}, the pseudo$\mbox{-}$circle family must create exactly two distinct envelopes.
	
	For case $(B)$, the proof is similar as the proof of $(1\mbox{-}\infty)$. The pseudo$\text{-}$circle $C_{(\bm{a}(t), \pm r(t))}$ creates uncountably many distinct envelopes if $\beta(t)\neq0$ is not dense in $I$.
	
	Therefore, it is impossible for the pseudo$\mbox{-}$circle family $C_{(\bm{a}(t), \pm r(t))}$ creating a unique envelope if $\epsilon_{\tilde{\bm{\nu}}}\cdot\epsilon_{\bm{\mu}}=-1$.\\
	\\	
	$(2\mbox{-}ii)$ Suppose that the set consisting of $t\in I$ satisfying $\beta(t)\neq0$ is not dense in $I$. By $(2\mbox{-}i ~case~(B))$, there are uncountably many distinct envelopes created by $C_{(\bm{a}(t), \pm r(t))}$. This contradicts the assumption that there are two distinct envelopes created by $C_{(\bm{a}(t), \pm r(t))}$. Thus, the set consisting of $t\in I$ satisfying $\beta(t)\neq0$ must be dense in $I$.
	
	Conversely, if the set consisting of $t\in I$ satisfying $\beta(t)\neq0$ is dense in $I$, we can prove the pseudo$\mbox{-}$circle family $C_{(\bm{a}(t), \pm r(t))}$ creates two distinct envelopes by the same way as the proof of  $(2\mbox{-}i ~case~(A))$.\\
	\\	
	$(2\mbox{-}\infty)$ The proof is similar to the proof of $(1\mbox{-}\infty)$.
\end{proof}
Compared with the assertion $(3)$ of Theorem $1$ in \cite{Wt}, the properties of envelopes created by pseudo-circle families exhibit variability and differences. For instance, from the $(2\mbox{-}i)$ of Theorem \ref{1.12}, it can be inferred that the pseudo$\mbox{-}$circle family $C_{(\bm{a}(t), \pm r(t))}$ cannot create a unique envelope if $\epsilon_{\tilde{\bm{\nu}}}\cdot\epsilon_{\bm{\mu}}=-1$.

\section{Examples\label{section3}}
\begin{example}\label{3.1}
	We consider the Example \ref{1.111} by applying Theorems \ref{1.14}, \ref{1.11} and \ref{1.12}. In Example \ref{1.111}, we define the timelike frontal $\bm{a}:I\rightarrow\mathbb{R}_1^2$ as $\bm{a}(t)=(t^3,\sqrt{1+t^6})$, $t\in\mathbb{R}$. Thus, the mapping $\bm{\nu}:\mathbb{R}\rightarrow S_1^1$ and $\bm{\mu}:\mathbb{R}\rightarrow H^1$ are obtained as $\bm{\nu}(t)=(t^3,\sqrt{1+t^6})$ and $\bm{\mu}(t)=(\sqrt{1+t^6},t^3)$, respectively. The radius function $r:\mathbb{R}\rightarrow\mathbb{R}^+$ is given by $r(t)=1$. Consider the pseudo$\mbox{-}$circle family $C_{(\bm{a}(t), r(t))}$, we have
	$$\frac{d r}{dt}(t)=0,$$
	and
	$$\beta(t)=-\langle \frac{d \bm{a}}{dt}(t),\bm{\mu}(t)\rangle=\frac{3t^2}{\sqrt{1+t^6}}.$$
	Since $\epsilon_{\tilde{\bm{\nu}}}\cdot\epsilon_{\bm{\mu}}=-1$, the creative condition is as follows:
	$$\frac{d r}{dt}(t)=\beta(t)\sinh\theta(t).$$
	Thus, take $\sinh\theta(t)=0$ and $\cosh\theta(t)=1$ such that $C_{(\bm{a}(t), r(t))}$ is a creative. By Theorem \ref{1.14}, the pseudo$\mbox{-}$circle family $C_{(\bm{a}(t), r(t))}$ creates an envelope $f:I\rightarrow\mathbb{R}^2_1$. At this case, the unit vector $\tilde{\bm{\nu}}\in S_1^1$ satisfying
	$$\tilde{\bm{\nu}}=\pm\bm{\nu}(t)=\pm(t^3,\sqrt{1+t^6}).$$
	Therefore, according to Theorem \ref{1.11}, the envelopes $f(t)$ created by $C_{(\bm{a}(t), r(t))}$ is parametrized as follows:
	\begin{align*}
		f(t)&=\bm{a}(t)+r(t)\tilde{\bm{\nu}}(t)\\
		&=(t^3,\sqrt{1+t^6})\pm(t^3,\sqrt{1+t^6})\\
		&=(t^3\pm t^3,\sqrt{1+t^6}\pm\sqrt{1+t^6}).	
	\end{align*}
	Finally, by the assertion $(2)$ of Theorem \ref{1.12}, the pseudo$\mbox{-}$circle family $C_{(\bm{a}(t), r(t))}$ creates exactly two envelopes.
\end{example}
\begin{example}
	Let $\bm{a}:I\rightarrow\mathbb{R}_1^2$ be a spacelike frontal given by $\bm{a}(t)=(\cosh(t),\sinh(t))$, where $t\in(-2,2)$. Thus, if we define the unit vector $\bm{\nu}(t)=(\cosh(t),\sinh(t))$,  $\bm{\nu}:(-2,2)\rightarrow H^1$ gives the Gauss mapping of $\bm{a}$. By definition, $\bm{\mu}(t)=(\sinh(t),\cosh(t))$ and we have $\beta(t)=\langle \frac{d \bm{a}}{dt}(t),\bm{\mu}(t)\rangle=1$. Furthermore, the radius function $r:I\rightarrow\mathbb{R}^+$ is defined by $r(t)=2-t$. Consider the pseudo$\mbox{-}$circle family $C_{(\bm{a}(t), r(t))}$, we have
	$$\frac{d r}{dt}(t)=-1.$$
	Since $\epsilon_{\tilde{\bm{\nu}}}\cdot\epsilon_{\bm{\mu}}=1$, the creative condition is $\frac{d r}{dt}(t)=-\beta(t)\cosh\theta(t)$. It follows from $C_{(\bm{a}(t), r(t))}$ is creative that $\cosh\theta(t)=1$ and $\sinh\theta(t)=0$. According to the assertion $(1)$ of Theorem \ref{1.11}, we obtain
	$$\tilde{\bm{\nu}}(t)=\bm{\mu}(t)=(\sinh(t),\cosh(t)).$$
	Thus, we have
	$$f(t)=\bm{a}(t)+r(t)\tilde{\bm{\nu}}(t)=(\cosh(t)+(2-t)\sinh(t), \sinh(t)+(2-t)\cosh(t)).$$
	Finally, by the assertion $(1)$ of Theorem \ref{1.12}, $f(t)$ is the unique envelope created by the pseudo$\mbox{-}$circle family $C_{(\bm{a}(t), r(t))}$, see Figure \ref{figure_example2}.
	\begin{figure}[h]
		\begin{center}
			\includegraphics[width=8cm]
			{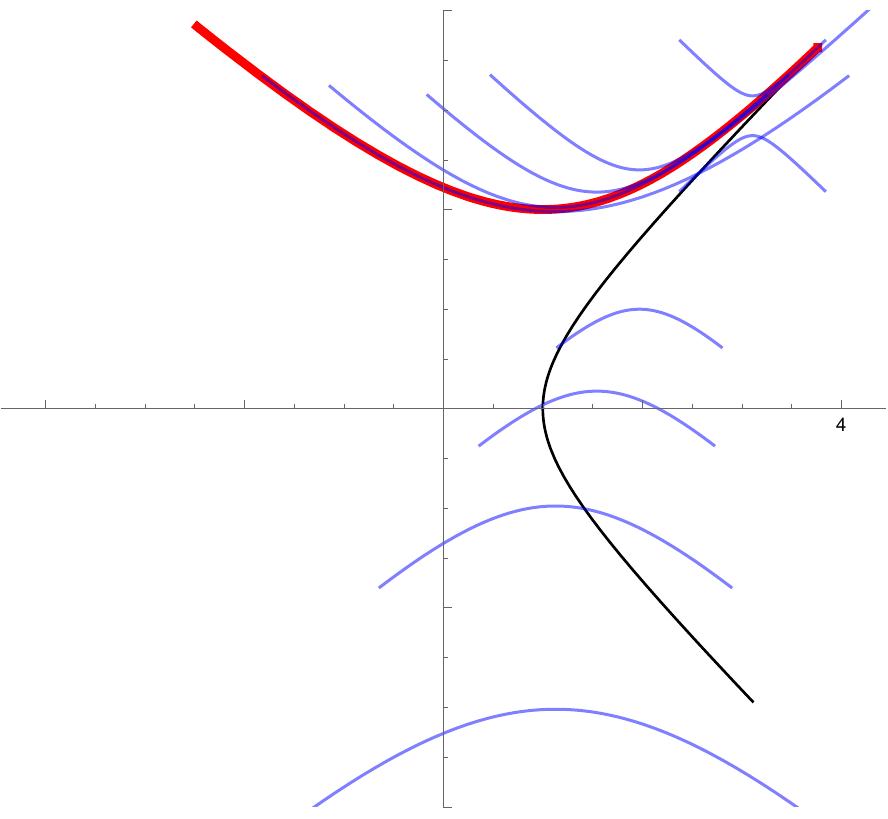}
			\caption{The pseudo$\mbox{-}$circle family $C_{(\bm{a}(t), r(t))}$
				and the candidates of its envelope.
			}
			\label{figure_example2}
		\end{center}
	\end{figure}
\end{example}
\begin{example}\label{2.11}
	We provide an example to demonstrate Proposition \ref{1.13}. Let $\bm{a}:I\rightarrow\mathbb{R}_1^2$ be a spacelike frontal defined by  $\bm{a}(t)=\big(\frac{t^2}{2},\frac{t^2}{2}+\frac{t^3}{3}\big)$ and the radius function $r:I\rightarrow\mathbb{R}^+$ defined by $r(t)=1$, where $I=(0,\infty)$. We consider the pseudo$\mbox{-}$circle family $C_{(\bm{a}(t), r(t))}$. Thus, we have $\bm{\nu}(t)=\frac{1}{\sqrt{t(t+2)}}(t+1,1)$ and $\bm{\mu}(t)=\frac{1}{\sqrt{t(t+2)}}(1,t+1)$. It follows  $\beta(t)=\langle \frac{d \bm{a}}{dt}(t),\bm{\mu}(t)\rangle=\frac{t^2(t+2)}{\sqrt{t(t+2)}}$. We want to find a mapping  $\tilde{\bm{\nu}}:I\rightarrow S^1_1$ satisfying
	$$\frac{d r}{dt}(t)=-\beta(t)\langle\tilde{\bm{\nu}}(t),\bm{\mu}(t)\rangle.$$ Namely, the mapping $\tilde{\bm{\nu}}:I\rightarrow S^1_1$ satisfying
	$$0=-\frac{t^2(t+2)}{\sqrt{t(t+2)}}\langle\tilde{\bm{\nu}}(t),\bm{\mu}(t)\rangle.$$
	Since $\epsilon_{\tilde{\bm{\nu}}}\cdot\epsilon_{\bm{\mu}}=1$, it follows that the creative conditions
	$$\frac{d r}{dt}(t)=-\beta(t)\cosh\theta(t)$$
	or
	$$\frac{d r}{dt}(t)=\beta(t)\cosh\theta(t)$$
	simply becomes
	$$0=-\frac{t^2(t+2)}{\sqrt{t(t+2)}}\cosh\theta(t)$$ or
	$$0=\frac{t^2(t+2)}{\sqrt{t(t+2)}}\cosh\theta(t).$$
	Since the function $\cosh\theta(t)\geq1$, the above equations do not hold for all $t\in I$. By Theorem \ref{1.14}, the pseudo$\mbox{-}$circle family $C_{(\bm{a}(t), r(t))}$ does not create an envelopes.
\end{example}
\begin{example}\label{2.12}
	In this example, we define $\bm{a}(t)=\big(\frac{t^2}{2},\frac{t^2}{2}+\frac{t^3}{3}\big)$ and $r(t)=1$, where $t\in(-2,0)$. The difference from Example \ref{2.11} is the parameter space, where $I$ is now $(-2,0)$ instead of $(0,\infty)$. We will continue to consider the pseudo$\mbox{-}$circle family $C_{(\bm{a}(t), r(t))}$ with the modified parameter space.
	
	We defined $\bm{\nu}(t)=\frac{1}{\sqrt{-t(t+2)}}(t+1,1)$ and $\bm{\mu}(t)=\frac{1}{\sqrt{-t(t+2)}}(1,t+1)$. It follows $\bm{a}:I\rightarrow\mathbb{R}^2_1$ is a timelike frontal and $\beta(t)=-\langle \frac{d \bm{a}}{dt}(t),\bm{\mu}(t)\rangle=-\frac{t^2(t+2)}{\sqrt{-t(t+2)}}$. Since $\epsilon_{\tilde{\bm{\nu}}}\cdot\epsilon_{\bm{\mu}}=-1$, it follows that the  equivalent creative condition is
	$$\frac{d r}{dt}(t)=\beta(t)\sinh\theta(t).$$
	It simply becomes
	$$0=-\frac{t^2(t+2)}{\sqrt{-t(t+2)}}\sinh\theta(t).$$
	Since the above equation holds for any $t\in I$ if $\sinh\theta(t)=0$, it follows that the pseudo$\mbox{-}$circle family $C_{(\bm{a}(t), r(t))}$ is creative. Therefore, according to Theorem 1, $C_{(\bm{a}(t), r(t))}$ creates an envelope. Additionally, according to Theorem \ref{1.11}, the expression of $\tilde{\bm{\nu}}(t)$ must be as follows:
	$$\tilde{\bm{\nu}}(t)=\sinh\theta(t)\bm{\mu}(t)\pm\cosh\theta(t)\bm{\nu}(t)=\pm\bm{\nu}(t).$$
	Therefore, an envelope $f(t)$ created by $C_{(a(t), r(t))}$ is parameterized as follows:
	\begin{align*}
		f(t)&=\bm{a}(t)+r(t)\tilde{\bm{\nu}}(t)\\
		&=\bigg(\frac{t^2}{2},\frac{t^2}{2}+\frac{t^3}{3}\bigg)\pm\frac{1}{\sqrt{-t(t+2)}}(t+1,1)\\
		&=\bigg(\frac{t^2}{2}\pm\frac{t+1}{\sqrt{-t(t+2)}},\frac{t^2}{2}+\frac{t^3}{3}\pm\frac{1}{\sqrt{-t(t+2)}}\bigg).\\	
	\end{align*}
	Finally, by the assertion $(2)$ of Theorem \ref{1.12}, the number of envelopes created by the pseudo-circle family $C_{(\bm{a}(t), r(t))}$ is exactly two, see Figure \ref{figure_example3}.
	\begin{figure}[h]
		\begin{center}
			\includegraphics[width=8cm]
			{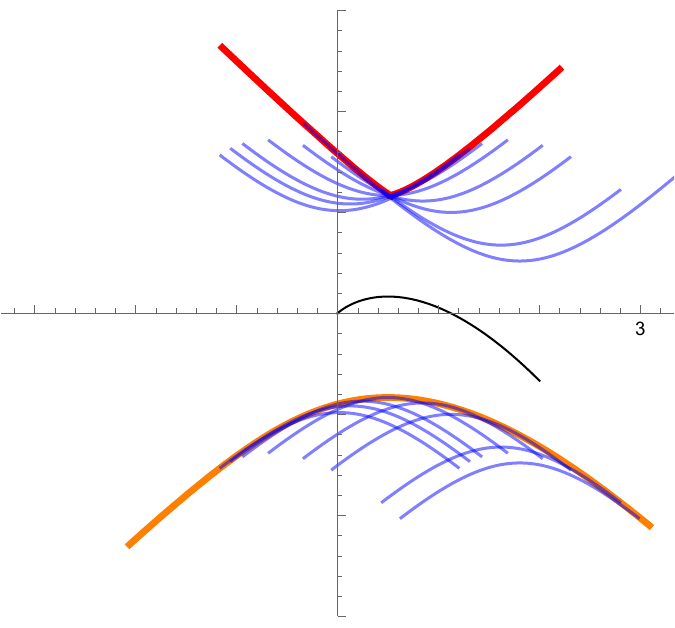}
			\caption{The pseudo$\mbox{-}$circle family $C_{(\bm{a}(t), r(t))}$
				and the candidates of its envelope.
			}
			\label{figure_example3}
		\end{center}
	\end{figure}
\end{example}
\begin{example}
	In this example, we set $\bm{a}(t)=(0,1)$ and $r(t)=1$, where $t\in\mathbb{R}$. Thus, every mapping $\bm{\nu}:\mathbb{R}\rightarrow H^1$ can be chosen as the Gauss mapping such that $\bm{a}$ is a frontal.  We consider the pseudo$\mbox{-}$circle family $C_{(\bm{a}(t), r(t))}$. We have $\beta(t)=\langle \frac{d \bm{a}}{dt}(t),\bm{\mu}(t)\rangle=0$. Since the radius function $r(t)$ is a constant function, the creative condition
	$$\frac{d r}{dt}(t)=-\beta(t)\langle\tilde{\bm{\nu}}(t),\bm{\mu}(t)\rangle$$ simply becomes
	$$0=0.$$
	Therefore, the creative condition is satisfied for any $\tilde{\bm{\nu}}:I\rightarrow S^1_1$. According to Theorem \ref{1.14}, the pseudo$\mbox{-}$circle family $C_{(\bm{a}(t), r(t))}$ creates an envelope. According to Theorem \ref{1.11}, the parametrization of the envelope is given by
	$$f(t)=\bm{a}(t)+r(t)\tilde{\bm{\nu}}(t)=(0,1)+\tilde{\bm{\nu}}(t).$$
	Finally, by Theorem \ref{1.12}, there are uncountably many distinct envelope created by $C_{(\bm{a}(t), r(t))}$.
\end{example}
\section{Alternative definitions\label{section4}}
In Definition \ref{1.15} of Section \ref{section2}, the definition of envelope created by a pseudo$\mbox{-}$circle family $C_{(\bm{a}(t), \pm r(t))}$ is provided. In \cite{Wt}, the set
comprising the images of envelopes of circle families defined in Definition $2$ is referred to as $E_2$ \textit{envelope}. Naturally, we call the envelope of pseudo-circle families defined in Definition \ref{1.111} an $E_2$ \textit{envelope}. Additionally, two alternative definitions, namely the $E_1$ \textit{envelope} and $\mathcal{D}$ \textit{envelope} for pesudo-circle families, are provided as follows.
\begin{definition}
	($E_1$ \textit{envelope}). Consider a frontal, denoted by $\bm{a}:I\rightarrow\mathbb{R}^2_1$ and let $r:I\rightarrow\mathbb{R}^+$ be a positive function. Consider a parameter $t_0$ in the interval $I$ and fix it. Suppose that
	$$\lim_{\varepsilon\rightarrow0}C_{(\bm{a}(t_0), \pm r(t_0))}\cap C_{(\bm{a}(t_0+\varepsilon), \pm r(t_0+\varepsilon))}$$
	is not an empty set. We denote this set as $I(t_0)$. Consider a point $\bm{e}_1(t_0)=(x(t_0),y(t_0))$ belonging to the set $I(t_0)$. Then, the set comprising the images of smooth mappings $\bm{e}_1:I\rightarrow\mathbb{R}^2_1$ is referred to as the $E_1$ \textit{envelope} created by the pseudo$\mbox{-}$circle family $C_{(\bm{a}(t), \pm r(t))}$ and it is denoted by $E_1$.
\end{definition}
\begin{definition}\label{4.1}
	($\mathcal{D}$ \textit{envelope}). Consider a frontal denoted by $\bm{a}:I\rightarrow\mathbb{R}^2_1$ and let $r:I\rightarrow\mathbb{R}^+$ be a positive function. Set
	$$F(x,y,t)=\|(x,y)-\bm{a}(t)\|^2\mp r(t)^2,$$
	then the set defined by the conditions
	$$\left\{(x, y)\in \mathbb{R}^2_1\, \left|\,
	\exists t \mbox{ such that }F(x,y,t)=\frac{\partial F}{\partial t}(x, y, t)=0
	\right.\right\}$$
	is referred to as $\mathcal{D}$ \textit{envelope} created by the pseudo$\mbox{-}$circle family $C_{(\bm{a}(t), \pm r(t))}$ and it is denoted by $\mathcal{D}$.
\end{definition}

We will study the relationship among $E_1$, $E_2$ and $\mathcal{D}$. The relationship between the $E_1$ \textit{envelope} and $E_2$ \textit{envelope} is as follows.
\begin{theorem}
	$E_1=E_2.$
\end{theorem}
\begin{proof}
	To prove that $E_1\subset E_2$, we consider a  parameter $t_0$ in the interval $I$ and a sequence ${t_i}$ in $I$ conversing to $t_0$, where $i=1,2,\ldots$. Let $(x(t_0),y(t_0))$ be a point in $E_1$. We assume that the point $(x(t_i),y(t_i))$ is an intersection point of two pseudo$\mbox{-}$circle family $C_{(\bm{a}(t_i), \pm r(t_i))}\cap C_{(\bm{a}(t_0), \pm r(t_0))}$ and satisfies the condition
	$$\lim_{t_i\rightarrow t_0}(x(t_i),y(t_i))=(x(t_0),y(t_0)).$$
	Then, we have the following equations:
	$$\|(x(t_i),y(t_i))-\bm{a}(t_i)\|^2=\pm(r(t_i))^2,$$
	$$\|(x(t_i),y(t_i))-\bm{a}(t_0)\|^2=\pm(r(t_0))^2.$$
	Set $\bm{a}(t_j)=(a_x(t_i),a_y(t_i))$ for any $j=0,1,2,\ldots.$, it follows that
	$$2\bigg(x(t_i)\big(a_x(t_i)-a_x(t_0)\big)-y(t_i)\big(a_y(t_i)-a_y(t_0)\big)\bigg)-a_x^2(t_i)+a_x^2(t_0)+a_y^2(t_i)-a_y^2(t_0)=\pm\big((r(t_i))^2-(r(t_0))^2\big).$$
	Since $\lim_{i\rightarrow\infty}t_i=t_0$ and $\lim_{t_i\rightarrow t_0}\big(x(t_i),y(t_i)\big)=\big(x(t_0),y(t_0)\big)$, we have following
	$$2\big(x(t_0)\frac{d a_x}{dt}(t_0)-y(t_0)\frac{d a_y}{dt}(t_0)\big)-2a_x(t_0)\frac{d a_x}{dt}(t_0)+2a_y(t_0)\frac{d a_y}{dt}(t_0)=\pm2r(t_0)\frac{d r}{dt}(t_0).$$
	Then
	$$\mp\frac{1}{r(t_0)}\bigg\langle\big(x(t_0)-a_x(t_0),y(t_0)-a_y(t_0)\big),\big(\frac{d a_x}{dt}(t_0),\frac{d a_y}{dt}(t_0)\big)\bigg\rangle=\frac{d r}{dt}(t_0).$$
	Since the vector
	$$\frac{1}{r(t_0)}\big(x(t_0)-a_x(t_0),y(t_0)-a_y(t_0)\big)=\frac{1}{r(t_0)}\big((x(t_0),y(t_0))-\bm{a}(t_0)\big)$$
	is a unit vector and $\big(\frac{d a_x}{dt}(t_0),\frac{d a_y}{dt}(t_0)\big)=\beta(t_0)\bm{\mu}(t_0)$. As a result, the condition for creative is satisfied at $t=t_0$. Consequently, according to Theorem \ref{1.14}, the point $(x(t_0),y(t_0))$ belongs to $E_2$.
	
	On the other hand, if the pseudo$\mbox{-}$circle family $C_{(\bm{a}(t),\pm r(t))}$ creates an $E_2$ envelope $f:I\rightarrow\mathbb{R}^2_1$, then by the assertion of Theorem \ref{1.11}, the representation of $f$ is as follows
	$$f(t)=\bm{a}(t)+r(t)\tilde{\bm{\nu}}(t).$$
	We define the straight line $L(P,\bm{\nu})$ for a point $P\in\mathbb{R}^2_1$ and a unit vector $\bm{\nu}(t):I\rightarrow H^1(\rm{or }~ S^1_1)$ as follows:
	$$L(P,\bm{\nu})=\{(x,y)\in\mathbb{R}^2_1~|~\big\langle((x,y)-P),\bm{\nu}\big\rangle=0\}.$$
	Since
	\begin{align*}
		&\big\langle \frac{d f}{dt}(t),\tilde{\bm{\nu}}(t)\big\rangle\\
=&\big\langle \frac{d \bm{a}}{dt}(t)+\frac{d r}{dt}(t)\tilde{\bm{\nu}}(t)+r(t)\frac{d \tilde{\bm{\nu}}}{dt}(t),\tilde{\bm{\nu}}(t)\big\rangle\\
=&\langle \beta(t)\bm{\mu}(t),\tilde{\bm{\nu}}(t)\rangle\pm \frac{d r}{dt}(t)\\
=&0,
	\end{align*}
	similar to Definition 1 in \cite{nishimura}, we call $f$ an $E_2$ envelope created by the straight line family
	$$\mathcal{L}(f,\tilde{\bm{\nu}})=\{L_{(f(t),\tilde{\bm{\nu}}(t))}\}_{t\in\mathbb{R}}.$$
	Let $t_0$ be a parameter in the interval $I$ and consider a sequence ${t_i}$ in $I$ conversing to $t_0$, where $i=1,2,\ldots$. By analysis similar to Theorem $2$ in \cite{Wt}, for sufficiently large $i\in\mathbb{N}$ there exists a point
	$$(x(t_i),y(t_i))\in L_{(f(t_0),\tilde{\bm{\nu}}(t_0))}\cap L_{(f(t_i),\tilde{\bm{\nu}}(t_i))}$$
	such that $\lim_{i\rightarrow\infty}(x(t_i),y(t_i))=f(t_0)$. This implies that for any sufficiently large $i\in\mathbb{N}$, there must exist a point
	$$(\tilde{x}(t_i),\tilde{y}(t_i))\in C_{(\bm{a}(t_0),\pm r(t_0))}\cap C_{(\bm{a}(t_i),\pm r(t_i))}$$
	such that $\lim_{i\rightarrow\infty}(\tilde{x}(t_i),\tilde{y}(t_i))=f(t_0)$.
	Therefore, we conclude that the point $f(t_0)\in\mathbb{R}^2_1$ belongs to $E_1$. Since $f$ is an arbitrary envelope created by $C_{(\bm{a}(t),\pm r(t))}$ and $t_0$ is an  arbitrary parameter in $I$, it follows that $E_2\subset E_1$.
\end{proof}
We prove the following theorem which asserts that $\mathcal{D}=E_2$ if and only if $\bm{a}(t)$ is regular, and $\mathcal{D}$ contains not only $E_2$ but also the pseudo-circle $C_{(\bm{a}(t_0), \pm r(t_0))}$ at a singular point $t_0$ of $\bm{a}(t)$ if $\bm{a}(t)$ is singular.
\begin{theorem}
	Let $\bm{a}:I\rightarrow\mathbb{R}^2_1$ be a frontal in the Minkowski plane and	$r:I\rightarrow\mathbb{R}^+$ a positive function.
	Suppose that pseudo-circle family $C_{(\bm{a}(t), \pm r(t))}$ creates the envelope $f$, then we have the following conclusions.\\
	$(1)$~$\mathcal{D}=f$ if $\bm{a}(t)$ is regular for any $t\in I$.\\
	$(2)$~$\mathcal{D}=f\cup C_{(\bm{a}(t_1), \pm r(t_1))}\cup\cdots\cup C_{(\bm{a}(t_k), \pm r(t_k))}$	if $t_1,\ldots,t_k$ are the singular points of $\bm{a}(t)$.
\end{theorem}
\begin{proof}
	Recalling that 	$\mathcal{D}$ is defined as
	$$\mathcal{D}=\left\{(x, y)\in \mathbb{R}^2_1\, \left|\,
	\exists t \mbox{ such that }F(x,y,t)=\frac{\partial F}{\partial t}(x, y, t)=0
	\right.\right\}.$$
	Let $(x_0,y_0)$ be a point of $\mathcal{D}$. Since $F(x,y,t)=\|(x,y)-\bm{a}(t)\|^2\mp r(t)^2$,  it
	follows that there exists a $t_0\in I $ such that the following $(a)$ and $(b)$ are satisfied.
	
	$(a)$~~$\langle(x_0,y_0)-\bm{a}(t_0),(x_0,y_0)-\bm{a}(t_0)\rangle\mp r(t_0)^2=0.$
	
	$(b)$~~$\frac{d\big(\langle(x_0,y_0)-\bm{a}(t_0),(x_0,y_0)-\bm{a}(t_0)\rangle\mp r(t_0)^2\big)}{dt}=0.$
\\	
The condition $(a)$ follows the existence of a unit vector field $\bm{\nu}_1(t):I\rightarrow S^1_1(\rm{or }~ H^1)$ satisfying
	$$(x_0,y_0)=\bm{a}(t_0)+r(t_0)\bm{\nu}_1(t_0).$$
	Condition $(b)$ derives the following equation
	$$\big\langle \frac{d \bm{a}}{dt}(t_0),(x_0,y_0)-\bm{a}(t_0)\big\rangle\pm r(t_0)\frac{d r}{dt}(t_0)=0.$$
	By substituting $\frac{d \bm{a}}{dt}(t_0)=\beta(t_0)\bm{\mu}(t_0)$, we obtain
	$$r(t_0)\big(\beta(t_0)\langle\bm{\mu}(t_0),\bm{\nu}_1(t_0)\rangle\pm \frac{d r}{dt}(t_0)\big)=0.$$
	Since $r(t)>0$ for any $t\in I$, we can conclude that
	$$\frac{d r}{dt}(t_0)=\mp\beta(t_0)\langle\bm{\mu}(t_0),\bm{\nu}_1(t_0)\rangle.$$
	
	On the other hand, since $C_{(\bm{a}(t), \pm r(t))}$ is a creative pseudo-circle family, there exists a smooth unit vector field $\tilde{\bm{\nu}}:I\rightarrow S^1_1(\rm{or }~ H^1)$ along $\bm{a}:I\rightarrow\mathbb{R}^2_1$ such that
	$$\frac{d r}{dt}(t)=\mp\beta(t)\langle\bm{\mu}(t),\tilde{\bm{\nu}}(t)\rangle$$
	for any $t\in I$. Assuming that the parameter $t_0\in I$ is a regular point of $\bm{a}(t)$, we have $\beta(t_0)\neq0$. Therefore, by the proof of Theorem \ref{1.14}, we conclude that
	$$(x_0,y_0)\in E_2.$$
	Suppose that the parameter $t_0\in I$ is a singular point of $\bm{a}(t)$, we have $\beta(t_0)=0$. Thus, for any unit vector $\bm{v}:I\rightarrow S^1_1(\rm{or }~ H^1)$, we have
	$$\frac{d r}{dt}(t_0)=\mp\beta(t_0)\langle\bm{\mu}(t_0),\bm{v}(t_0)\rangle.$$
	Hence, at a singular point $t_0\in I$ of $\bm{a}(t)$, we can select any unit vector $\bm{v}(t)\in S^1_1(\rm{or }~ H^1)$ as the unit vector $\bm{\nu}_1(t_0)$. Therefore, it follows
	$$\mathcal{D}_0=C_{(\bm{a}(t_0), \pm r(t_0))},$$
	where $\mathcal{D}_0=\left\{(x, y)\in \mathbb{R}^2_1\, \left|\,
	F(x,y,t_0)=\frac{\partial F}{\partial t}(x, y, t_0)=0
	\right.\right\}.$
\end{proof}
\section*{Acknowledgement}
The first author and the fourth author would like to thank Professor Takashi Nishimura for his help and guidance during their stay at Yokohama National University.
\section*{Funding}
The first author was supported by the National Natural Science Foundation of China
(Grant No. 12001079).

\end{document}